May 10, 2005

On Powers of a Hypergeometric Function

by
Michael Milgram[*]

Consulting Physicist, Geometrics Unlimited, Ltd.
Box 1484, Deep River, Ont. Canada. K0J 1P0

**Abstract**: A new expansion for integral powers of the hypergeometric function corresponding to a special case of the incomplete beta function is summarized, and consequences, including two new sums involving digamma (psi) functions are presented.

In another publication[1] (I), a new series expansion for the *arctangent* function was developed, and consequences noted. The motivation for that work, presented elsewhere[2], was a need to obtain closed form representations for the Fourier coefficients of powers of *arctangent*, appearing in applications. Fundamental to that analysis was the necessity to analytically sum various slowly converging power series involving digamma and polygamma functions; that work has also been reported elsewhere[3]. After these works were completed, I belatedly noticed that the methods described could be generalized to a special case of the incomplete beta function[†] $_2F_1(1,\alpha;\alpha+1;x)$, where *arctangent* corresponds to $\alpha = \frac{1}{2}$. This note is a record of the new, generalized results.

Following I, we have

$$\sum_{n=0}^{k} \frac{1}{(k-n+\alpha)} \sum_{m=0}^{n} \frac{1}{(m+\alpha)} = 2 \sum_{n=0}^{k} \frac{1}{(n+2\alpha)} \sum_{m=0}^{n} \frac{1}{(m+\alpha)} \qquad (1)$$

The proof, by induction, follows exactly the same steps as in I (Theorem 1) with the generalization $\frac{1}{2} \to \alpha$, leading to the following generalization of Theorem 2 of that paper

$$\sum_{m_n=0}^{k} \frac{1}{(k-m_n+\alpha)} \sum_{m_{n-1}=0}^{m_n} \frac{1}{(m_{n-1}+(n-1)\alpha)} \cdots \sum_{m_2=0}^{m_3} \frac{1}{(m_2+2\alpha)} \sum_{m_1=0}^{m_2} \frac{1}{(m_1+\alpha)}$$
$$= n \sum_{m_n=0}^{k} \frac{1}{(m_n+n\alpha)} \sum_{m_{n-1}=0}^{m_n} \frac{1}{(m_{n-1}+(n-1)\alpha)} \cdots \sum_{m_2=0}^{m_3} \frac{1}{(m_2+2\alpha)} \sum_{m_1=0}^{m_2} \frac{1}{(m_1+\alpha)} \ . \qquad (2)$$

To obtain this result, it is necessary to introduce recursively defined coefficients

$$t_m^n(\alpha) = \sum_{k=0}^{m} \frac{1}{k+n\alpha} t_k^{n-1}(\alpha) \qquad (3)$$

with

---
[*] mike@geometrics-unlimited.com
[†] Lerch's Φ function

$$t_0^n(\alpha) = \frac{1}{n!\alpha^n}$$
$$t_m^0(\alpha) = 1$$
$$t_k^1(\alpha) = \psi(k+1+\alpha) - \psi(\alpha) ,$$
(4)

from which the inductive proof follows exactly the same steps as in I.

Using partial fraction decomposition and with reference to (2), obtain

$$\sum_{m=0}^{k} \frac{1}{(k-m+\alpha)} \frac{1}{(m+n\alpha)} t_m^{n-1}(\alpha) = \frac{n+1}{k+(n+1)\alpha} \sum_{m=0}^{k} \frac{1}{(m+n\alpha)} t_m^{n-1}(\alpha)$$
$$= \frac{n+1}{k+(n+1)\alpha} t_k^n(\alpha).$$
(5)

Now consider a special case of the incomplete beta function:

$$({}_2F_1(1,\alpha;\alpha+1;x))^n = \left(\alpha \sum_{m=0}^{\infty} \frac{x^m}{\alpha+m}\right)^n = \alpha^n \sum_{m_n=0}^{\infty} x^{m_n} \sum_{m_{n-1}=0}^{m_n} \frac{1}{\alpha+m_n-m_{n-1}} \cdots$$
$$\times \sum_{m_2=0}^{m_3} \frac{1}{\alpha+m_3-m_2} \sum_{m_1=0}^{m_2} \frac{1}{(\alpha+m_1)(\alpha+m_2-m_1)}$$
(6)

using well-known (**4** Section 4.1) indicial shift operations for multiple sums. Apply (5) repeatedly from the right increasing "*n*" by one at each application, to find

$$({}_2F_1(1,\alpha;\alpha+1;x))^n = \alpha^n n! \sum_{m=0}^{\infty} \frac{x^m}{m+n\alpha} t_m^{n-1}(\alpha) \qquad \alpha \neq -m_n/n \ ; \ |x| \leq 1$$
$$= \alpha^n n! \sum_{0 \leq m_1 \leq m_2 \leq \ldots \leq m_n} \frac{x^{m_n}}{(m_n+n\alpha)(m_{n-1}+(n-1)\alpha)\cdots(m_1+\alpha)}$$
(7)

Lacking the motivation to evaluate specific integrals[2], (7) still has useful consequences. Consider the case *n=2*, which, after some re-ordering gives[5]

$$({}_2F_1(1,\alpha;\alpha+1;x))^2$$
$$= \alpha \sum_{l=0}^{\infty} \frac{1}{(l+\alpha)(l+\alpha+1)} {}_3F_2(2, 2\alpha, \alpha+1+l; 2\alpha+1, \alpha+2+l; x).$$
(8)

Now, set $x = -1$ to eventually obtain a new result [3]

$$\sum_{l=0}^{\infty}\frac{\psi(\frac{\alpha+l}{2})-\psi(\frac{\alpha+1+l}{2})}{(l+1-\alpha)} = \tfrac{1}{4}\left[\left[\psi(\tfrac{1+\alpha}{2})-\psi(\tfrac{\alpha}{2})\right]^2 + \left[\psi(\alpha)-\psi(\alpha-\tfrac{1}{2})\right]\left[\psi(1-\alpha)-\psi(\alpha)\right]\right] \quad (9)$$

which can be rewritten as

$$\sum_{l=0}^{\infty}\frac{\psi(q+l)-\psi(q+l+\tfrac{1}{2})}{(l+1-q)(l+\tfrac{1}{2}-q)} = \left[\psi(\tfrac{1}{2}+q)-\psi(q)\right]^2 \\ + 4\pi\cot(2\pi q)\left[\psi(2q)-\psi(2q-\tfrac{1}{2})\right] - \frac{\pi\cot(q\pi)}{(q-\tfrac{1}{2})} \quad (10)$$

by considering the even and odd terms of (9) individually. Also, by applying a linear transformation of variables ($x \to \frac{x}{(1-x)}$) to the left-hand side of (8) and a simplification[6] (Eq. 7.4.1.5) to the right-hand side, find

$$\left({}_2F_1(1,1;\alpha+1;x)\right)^2 = \frac{\alpha}{2\alpha-1}(\psi(\alpha)-\psi(1-\alpha))\,{}_2F_1(1,2;2\alpha+1;x) \\ - 2\alpha^2\sum_{l=0}^{\infty}\frac{1}{(l+\alpha)(l+1+\alpha)(l+1-\alpha)}\,{}_2F_1(1,2;\alpha+2+l;x) \quad (11)$$

Set $x = \tfrac{1}{2}$ and with reference to [**6**, Eq. 7.3.7(18)] eventually discover a new sum

$$\sum_{l=0}^{\infty}\frac{(-1)^l\psi(l+q)}{(l+2q-1)} = \frac{\psi(q)}{2}\left[\psi(q)-\psi(q-\tfrac{1}{2})\right] - \tfrac{1}{8}\left[\psi(\tfrac{1+q}{2})-\psi(\tfrac{q}{2})\right]^2 \quad (12)$$

reducing to a known result [**4**, Eq. 55.2.3] when q=1.

Generalizations to n>2 can be obtained using the methods set forth in the Appendix to ref. 2, but no simple (i.e. succinct) results were found.

References

[1] Milgram, M. "A Series Expansion for Integral Powers of Arctangent", Integral Transforms and Special Functions (to be published 2006); also available as http://www.arXiv.org/abs/math.CA/040633 (2004a).

[2] Milgram, M. "Identification and Properties of the Fundamental Expansion Functions for Neutron Transport in an Infinite Homogeneous Medium", Annals of Nuclear Energy (to be published, 2005).

[3] Milgram, M. "On Some Sums of Digamma and Polygamma Functions", submitted for publication and available as http://www.arXiv.org/abs/math.CA/040634 (2004b).

[4] Hansen, E.R, "*A Table of Series and Products*", Prentice-Hall (1975).

---